\documentclass{article}
\usepackage{amssymb}
\begin{document}

\newtheorem{thm}{Theorem}
\newtheorem{prop}{Proposition}
\newtheorem{lem}{Lemma}
\newtheorem{cor}{Corollary}
\newcommand{\bR}{{\bf R}}
\newcommand{\db}{\overline{\partial}}
\newcommand{\cL}{{\cal L}}
\newcommand{\cH}{{\cal H}}
\newcommand{\diverg}{{\rm div}}
\newcommand{\Ka}{{\rm Kahler}}
\newcommand{\cD}{{\cal D}}
\newcommand{\bC}{{\bf C}}
\newcommand{\ut}{\underline{t}}
\newcommand{\up}{\underline{p}}
\newcommand{\Vol}{{\rm Vol}}
\newcommand{\tk}{\tilde{k}}
%+Title
\title{Calabi-Yau metrics on Kummer surfaces as a model gluing problem}
\author{S. K. Donaldson}
\date{\today}
\maketitle
%-Title

%+Abstract
%\begin{abstract}
%    There is abstract text that you should replace with your own. 
%\end{abstract}
%-Abstract

%+Contents
%\tableofcontents
%-Contents

\section{Introduction}
This is an entirely expository article. The background is as follows.
\begin{itemize}
\item In 1976, Yau proved the Calabi conjecture, establishing the existence of Kahler metrics with vanishing Ricci curvature on many compact complex manifolds. The simplest examples, in complex dimension $2$, are K3 surfaces, and in particular the resolutions of {\it Kummer surfaces}.
\item An important  development in differential geometry over the past 30 years has been the use of \lq\lq gluing constructions'', creating  solutions to geometric problems, in some asymptotic regime, from appropriate building blocks. Many of the ideas,and applications, are due to Taubes---for example in the case of Yang-Mills instantons\cite{kn:T}. These constructions have been applied to a host of different  problems. A few examples are:   metrics of exceptional holonomy (Joyce \cite{kn:J}), holomorphic curves
 (Floer \cite{kn:F}, and others) constant mean curvature surfaces (Kapouleas \cite{kn:K}). A notable feature of these developments is that, while the  geometric contexts vary greatly,  many of the  analytical issues are the same.

\end{itemize}

The purpose of this article is to explain  a gluing construction for some Calabi-Yau metrics on K3 surfaces, resolutions of Kummer surfaces. This idea is by no means new and several such constructions are already available in the literature.
The first was done by  Topiwala \cite{kn:Top}. This used twistor theory, rather than PDE methods, but the underlying idea is the same. In the PDE framework, there are general theories of Kovalev and Singer \cite{kn:KS} and Arezzo and Pacard \cite{kn:AP} which can be used to treat the problem. To fit into Kovalev and Singer's set-up, one can exploit the fact that the desired metrics are \lq\lq self-dual'', and to fit into Arezzo and Pacard's that that they have constant scalar curvature. Yet another approach would be to take a product with a trivial flat $3$-torus, and then  fit into the higher-dimensional set-up of Joyce. On the other hand these papers (of Kovalev-Singer, Arezzo-Pacard  and Joyce) are all quite long, use a fair amount of machinery and are really aimed at  more difficult problems. So it seems worthwhile to give a short and elementary treatment of this simple case; partly for the interest of the result and partly as a model gluing problem where the  techniques can be illustrated in a relatively simple differential geometric setting. Despite the large body of work on these gluing techniques  in the literature, there are many interesting problems which have yet to be tackled and it is possible that the approach we take here may be useful in this way.

From a technical point of view the argument we give here probably does not  differ in any significant fundamental way from those of Kovalev-Singer and Arezzo-Pacard. The main difference from the latter is that we use the theory of manifolds with \lq\lq tubular ends'', rather than weighted function spaces on asymptotically Euclidean manifolds.

\

One should emphasise that these gluing methods only prove special cases of Yau's result---in a small asymptotic regime in the appropriate moduli space. On the other hand they do have he merit of giving an almost-explicit description of the metrics. 

\section{Technical background}
\subsection{Analysis}
We want to apply  the theory of translation-invariant elliptic operators on cylinders $M\times \bR$, where $M$ is a compact  Riemannian $(n-1)$-manifold. For brevity we just consider the case we will need which is the operator
$\Delta+1$, where $\Delta$ is the standard Laplace operator of the product metric  and we use the sign convention that $\Delta$ is a positive operator. For $p>1$ and integers $k\geq 0$, define Sobolev spaces $L^{p}_{k}$ on $M\times \bR$ by taking the completion of the smooth compactly supported functions under the usual norm. Then we have 
\begin{prop}
For any $p,k$ the map $\Delta+1:L^{p}_{k}\rightarrow L^{p}_{k}$ is an isomorphism.
\end{prop}

For our purposes we can take $p=2$. The proof is straightforward, assuming standard results about the compact manifold $M$. Given a smooth function $\rho$ of compact support we want to solve the equation $(\Delta+1)f=\rho$.  We can do this by separation of variables, writing $\rho=\sum \rho_{\lambda}(t)\phi_{\lambda}$ where $\phi_{\lambda}$ is an orthonormal basis of eigenfunctions of the Laplacian
$\Delta_{M}$ on $M$---so $(\Delta+1)\phi_{\lambda}= (1+\lambda)\phi_{\lambda}$. We seek a solution $f=\sum f_{\lambda}(t) \phi_{\lambda}$, so we need to solve the ODE's 
$$  -\frac{d^{2} f_{\lambda}}{dt^{2}}+(1+\lambda) f_{\lambda}=  \rho_{\lambda},$$
which can be done by standard elementary arguments. The solutions have exponential decay and integration-by-parts is valid, so that
$$   \int_{-\infty}^{\infty} \left( \frac{df_{\lambda}}{dt} \right)^{2}+ (1+\lambda) f_{\lambda}^{2}\ dt = \int_{-\infty}^{\infty} f_{\lambda} \rho_{\lambda}\  dt. $$
 Then the Cauchy-Schwartz inequality  implies that 
 a $$\int_{-\infty}^{\infty} f_{\lambda}^{2} \ dt \leq \int_{-\infty}^{\infty} \rho_{\lambda}^{2}\ dt, $$
and summing over $\lambda$ we see that the $L^{2}$ norm of the solution $f$ is bounded by that of $\rho$. Repeated integration by parts shows that for any $k$ there is a constant  $C_{k}$ such that we have $\Vert f\Vert_{L^{2}_{k+2}}\leq C_{k} \Vert \rho\Vert_{L^{2}_{k}}$ and the statement of the proposition (for $p=2$) is an easy consequence. (The integration-by-parts argument is made simpler if one uses the fact that on the compact manifold $M$ the $L^{2}_{k}$ norm is equivalent to
$$   \Vert g\Vert_{(k)}= \sum (\lambda+1)^{k} g_{\lambda}^{2} = \langle g,(\Delta+1)^{k} g\rangle, $$
for a function $g=\sum g_{\lambda}\phi_{\lambda}$.)

\

With the particular operator $\Delta+1$ the statement of the Proposition holds for very general class of manifolds, and can be proved in different ways. The advantage of the separation of variables approach above is that it extends easily to other elliptic operators on cylinders.

\

To tackle nonlinear problems we need Sobolev embedding theorems. These are easy to state.

\begin{prop}
If $k>l, k-n/p>l-n/q$ and $p<q$ then there is a continuous embedding $L^{p}_{k}\subset L^{q}_{l}$. If $k-n/p>0$ then there is a continuous embedding $L^{p}_{k}\subset C^{0}$. 
\end{prop}

Again, the proofs are not difficult, assuming facts about compact manifolds. Let us just consider the cases which will suffice in our application, when $n=4$. Then we want to establish  embeddings $L^{2}_{1}\subset L^{4}$ and $L^{2}_{3}\subset C^{0}$. For the first we use the fact that for functions $f$ on a \lq\lq band'' $M\times [0,1]$ we have an inequality
$$  \Vert f \Vert_{L^{4}}\leq C \Vert f \Vert_{L^{2}_{1}}. $$
(This follows from the usual theory for compact manifolds by considering the \lq\lq double'' of the band, i.e. $M\times S^{1}$.) Now decompose the cylinder $M\times \bR$ into a union of copies $\Omega_{n}=M\times [n,n+1]$ of the band.
If  $f$ is a function on $M\times \bR$ we  get
$$  \Vert f\Vert_{L^{4}}^{4} =\sum \int_{\Omega_{n}} f^{4} \leq C^{4} \sum \left( \int_{\Omega_{n}} \vert \nabla f \vert^{2}+f^{2}\right)^{2} \leq C^{4} \left( \int_{M\times \bR} \vert \nabla f\vert^{2}+ f^{2}\right)^{2}, $$
using the simple fact that for any $a_{n}\geq 0$ we have 
$$  \sum a_{n}^{2} \leq \left( \sum a_{n}\right)^{2}. $$

\

The inclusion $L^{2}_{3}\subset C^{0}$ is even easier--we simply multiply by a standard cut-off function supported in a band. A consequence of these two embeddings is that we have a bounded multiplication map
$  L^{2}_{3}\times L^{2}_{3} \rightarrow L^{2}_{3}. $

\

\ 

Now we move on to consider a Riemannian manifold $X$ with cylindrical ends, so the complement of a compact subset of $X$ is isometric to a finite disjoint union of half-cylinders.
$M_{i}\times (0,\infty)$. We consider an operator $\Box$ on $X$ of the form
$\Delta_{X}+ V$ where $V$ is a smooth function, equal to $1$ on each of the ends. We write $\cH_{\Box}$ for the set of functions $f$ in $L^{2}$ with $\Box f=0$. 
\begin{prop}\begin{enumerate}
\item  $\cH_{\Box}\subset L^{p}_{k}$ for all $p,k$.
 \item For any $p,k$ the operator $\Box:L^{p}_{k+2}\rightarrow L^{p}_{k}$ is Fredholm with kernel $\cH_{\Box}$ and image the orthogonal complement
 (in the $L^{2}$ sense) of $\cH_{\Box}$. 
\end{enumerate}
\end{prop}
 In fact  functions in $\cH_{\Box}$ have exponential decay, along with all their derivatives,  on the ends of the manifold. Usually one does not encounter manifolds with exactly cylindrical ends but rather ends which are asymptotic to cylinders (as Riemannian manifolds). The extension to this case is completely straightforward.

  \

  Now suppose we have a pair $X_{1}, X_{2}$ of such Riemannian manifolds with tubular ends. For simplicity of language, suppose that each has just one end and that the \lq\lq cross-section'' is the same compact manifold $M$. Given $T>0$ we form a compact manifold $X_{1}\sharp_{T}X_{2}$ by gluing the hypersurface corresponding to $M\times \{T\}$ in the end of $X_{1}$ to that in the end of $X_{2}$, in the obvious way. The result is a  Riemannian manifold which contains an isometric copy of $M\times (-T,T)$. Now suppose we have functions $V_{1}, V_{2}$ on $X_{1}, X_{2}$,as above. Then we get a function $V$ and an operator $\Box$ on $X_{1}\sharp_{T}V_{2}$ in the obvious way. (We use the same symbol $\Box$ to denote the operators on any of the manifolds involved.) The basic fact is

\begin{prop}
Suppose that $\Box$ is invertible on each of $X_{1}, X_{2}$. Then for any $p,k$ there is a constant $C_{p,k} $ and a $T_{0}$ such that if $T\geq T_{0}$ there is a right inverse $P$ to $\Box$ on $X_{1}\sharp_{T}X_{2}$ and $$  \Vert P \rho \Vert_{L^{p}_{k+2}}\leq C_{p,k} \Vert \rho\Vert_{L^{p}_{k}}.$$
\end{prop}

\

The crucial point here is that $C_{p,k}$ does not depend on $T$, once $T$ is sufficiently large.

  \

  The proof of this Proposition is simple. We fix a  partition of unity $\gamma_{1}+\gamma_{2}=1$ on $X_{1}\sharp_{T}X_{2}$with $\nabla \gamma_{i}$ supported in a standard band of width $1$ in the \lq\lq middle'' of the cylindrical region. Then we choose function $\beta_{1},\beta_{2}$ so that $\beta_{i}=1$ on the support of $\gamma_{i}$ but $\beta_{i}$ is supported in the region which can be considered, by an obvious stretch of langauge, as being contained in $X_{i}$. We choose $\beta_{i}$ so that $\nabla \beta_{i}$
is $O(T^{-1})$ and similarly for higher derivatives.  Let $P_{i}$ be the inverse  to $\Box$ over $X_{i}$ and set 
 $$  P_{0}\rho =\beta_{1} P_{1}(\gamma_{1} \rho) +\beta_{2}P_{2}(\gamma_{2}\rho), $$
where again we stretch notation to move  functions between $X_{i}$ and $X_{1}\sharp_{T}X_{2}$. Then $$\Box P_{0}\rho= \rho+ \sum_{i} 2 \nabla \beta_{i} \nabla P_{i}(\gamma_{i} \rho)+ \Delta \beta_{i} P_{i}(\gamma_{i} \rho), $$
and $$\Vert \Box P_{0}\rho-\rho \Vert_{L^{p}_{k+2}}\leq C T^{-1} \Vert \rho \Vert_{L^{p}_{k}},$$
so when $T$ is large enough we get a genuine right inverse $P_{0}\circ(\Box P_{0}-1)^{-1}$ and the estimate of the operator norm of $P$ is immediate.

\

The Sobolev embedding theorems on the infinite cylinder imply corresponding statements on $X_{1}\sharp_{T}X_{2}$, with constants independent of $T$.

\

\subsection{Geometry}

We recall some very standard facts about Kahler geometry,  the Kummer construction. and the Eguchi-Hanson metric.

Let $Z$ be a complex manifold of complex dimension $2$.  Giving a  Hermitian metric on $Z$ is the same as giving a positive form of type $(1,1)$. The metric is Kahler if this form is closed. Write
$\cD$ for the operator $2i\db\partial$ mapping (real) functions to (real) forms of type $(1,1)$.  If $\omega$ is a Kahler form the Laplacian of the metric is given by 
$$  \Delta_{\omega} f = (\cD f \wedge \omega)/\omega^{2}, $$
where \lq\lq division'' by the volume form $\omega^{2}$ has the obvious meaning. Suppose that $\chi$ is a nowhere-vanishing holomorphic $2$-form on $Z$. A Kahler metric is Calabi-Yau (i.e. Ricci-flat) if $\omega^{2}= \lambda \chi\wedge \overline{\chi}$, for some $\lambda>0$. If $\omega_{0}$ is one Kahler form and $\phi$ is a function then $\omega_{\phi}=\omega_{0}+\cD\phi$ is Kahler, provided it is positive (and positivity is an open condition). So we want to solve the Calabi-Yau equation
$$   (\omega_{0}+\cD\phi)^{2}= \lambda \chi\wedge \overline{\chi}, $$
with the side condition that $\omega_{0}+\cD\phi>0$.

\

Now we turn to the Kummer construction. Let $T^{4}=\bC^{2}/\Lambda$ be a complex torus. The map $z\mapsto-z$ on $\bC^{2}$ induces an involution of $T^{4}$ with $2^{4}=16$ fixed points. The quotient $\overline{X}$ is an orbifold with 16 singular points, each modelled on the quotient of $\bC^{2}$ by $\pm 1$. We write $X$ for the complement of the singular points in $\overline{X}$. The constant holomorphic $2$-form $dz_{1} dz_{2}$ is preserved by by the involution and so descends to a holomorphic form on  $X$.

Consider the map $(z_{1}, z_{2})\mapsto (z_{1}^{2}, z_{1} z_{2}, z_{2}^{2})\in \bC^{3}$. This induces a bijection between $\bC^{2}/\pm 1$ and the singular affine quadric in $\bC^{3}$ defined by the equation $v^{2}=uw$. We blow-up the origin in $\bC^{3}$ and take the proper transform of this affine surface in the blow-up. The result is a smooth surface $Y$, the resolution of this singularity. However all we really need to know is that $Y$ is a complex surface which, outside a compact set $K\subset Y$ is identified with the quotient  $(\bC^{2}\setminus
B^{4})/\pm 1$ and that the holomorphic form on this quotient extends to a nowhere-vanishing form on $Y$. Making this construction at each of the 16 singular points of $X$ we get a compact complex surface $Z$, with a nowhere vanishing holomorphic form. 

\

Now our gluing problem will be to find a Calabi-Yau metric on $Z$ starting with standard building blocks: metrics  on $X$ and $Y$. (Of course really we have 16 copies of $Y$.) The metric, $\omega_{X}$,  on $X$ that we need is just the flat one, but we also need a Calabi-Yau metric on $Y$. This is the {\it Eguchi-Hanson} metric, which we will now recall.

\

Go back to $\bC^{2}$ and write $\rho=r^{2}= \vert z_{1}\vert^{2} + \vert z_{2}\vert^{2}$.
Consider a Kahler metric of the form $\cD\psi$, where $\psi=F(\rho)$.
The Calabi-Yau equation becomes
$$  \det\left( \begin{array}{cc}F'+ \vert z_{1}\vert^{2} F'' & F'' z_{1} \overline{z}_{2}\\
F'' z_{2} \overline{z_{1}}& F'+ \vert z_{1}\vert^{2} F''\end{array} \right) = 1, $$
which is $ (F')^{2}+ \rho F' F''=1$. The solution $F(\rho)=\rho$ corresponds to the standard Euclidean metric $\Omega$. Up to re-scalings there is just one other solution which we can take to be given by
$$ F'(\rho) = \sqrt{ 1+\rho^{-2}}.$$
There is no need to integrate this explicitly, all we need is that, choosing the constant of integration suitably) we have
$ F(\rho)= \rho + G(\rho) $ where, for $\rho>1$
$G(\rho)$ has a convergent expansion $a_{1} \rho^{-1}+ a_{2} \rho^{-2}+\dots $
So we get a Calabi-Yau metric $ \Omega+ \cD G$ on $\bC^{2}\setminus \{0\}$ where $G=  a_{1} r^{-2}+a_{2} r^{-4}+ \dots$ for $r>1$.  
This metric has a singularity at the origin  but one can check that when we pass to the quotient and its resolution $Y$ we get a smooth Calabi-Yau metric $\omega_{Y}$. Choose a positive function $r_{Y}$ on $Y$ which is equal to $r=\sqrt{\vert z_{1}\vert^{2}+\vert z_{2}^{2}}$ on $Y\setminus K$. 

\

To set the scene for the gluing problem, fix a cut-off function $\beta$ on $\bR$, with $\beta(s)=0$ for $s\leq 1/2$ and $\beta(s)=1$ for $s\geq 1$. Introduce a (large) parameter $R$ and define a function $\beta_{R}$ on $Y$
by $\beta_{R}=  \beta( R^{-1/2} r_{Y})  $. Put
$$\omega_{R,Y}= \omega_{Y}- \cD(\beta_{R} G).  $$

Then, by construction, $\omega_{R,Y}$ equals the Eguchi-Hanson metric $\omega_{Y}$
when  $r_{Y}\leq \sqrt{R}/2$ and equals the flat metric $\Omega$ when $r_{Y}\geq \sqrt{R}$. The derivatives of $G$
satisfy
$$  \vert \nabla^{j} G \vert= O(r_{Y}^{-2-j}). $$
(Here we measure the size of derivatives with respect to the Euclidean metric.)
So on the annulus $\sqrt{R}/2\leq r_{Y}\leq \sqrt{R}$ we have $\vert \nabla^{j}(G)\vert=O(R^{-1-j/2})$.
The derivatives of $\beta_{R}$ satisfy (by scaling)
$$  \vert \nabla^{k} \beta_{R} \vert = O(R^{-k/2}),$$
so any product $\nabla^{j}\beta_{R}\nabla^{k} G$ is $O(R^{-1-(j+k)/2})$. Since $\cD(\beta_{R} G)$ is a sum of such products with $j+k=2$ we see that
$$  \vert \cD(\beta_{R} G)\vert = O(R^{-2}). $$
It follows, first,  that $\omega_{R,Y}$ is positive (for large enough $R$) so it is a Kahler metric. Second, we can write
$$  \omega_{R,Y}^{2} = (1+\eta)^{-1} \omega_{Y}^{2}$$ where  
$\eta$ is supported on this annulus and $\vert \eta\vert$ is $O(R^{-2})$. 

\

Now scale the metric  $\omega_{R,Y}$ by a factor  $R^{-2}$ ({\it i.e.} we scale lengths by a factor $R^{-1}$).  The sphere $r_{Y}=\sqrt{R}$ in  $Y$ is then isometric to a small sphere of radius $R^{-1/2}$ about each singular point in $X$.  Take 16 copies of $Y$, cut out 16 of these balls about the singular points, and glue in the corresponding region in the copies of $Y$. The result is  a Kahler metric $\omega_{0}$ on the complex manifold $Z$, which depends on the parameter $R$. (This parameter can be described more invariantly in terms of the  Kahler class. ) Our task is to deform this metric---the \lq\lq approximate solution''---to a genuine  Calabi-Yau metric on $Z$, once the parameter $R$ is sufficiently large.

\section{The gluing argument}
\subsection{Set-up} 

We want to treat our problem using the \lq\lq cylindrical ends'' theory and to do this we make a conformal change. The basic point is that $\bC^{2}\setminus \{0\}$ is conformally equivalent to the cylinder $S^{3}\times \bR$. However the metric on the cylinder is {\it not} Kahler. So consider in general a Kahler metric
$\omega$ and a positive real function $h$ on a complex surface and the conformally equivalent metric $\Theta= h^{-2} \omega$. WE write $d\mu$ for the volume form of the metric $\Theta$. Let $Q$ be the differential operator
$$  Q f= h \cD (h^{-1} f)  . $$
Notice that $Q$ is not changed if we multiply $h$ by a constant. Set
$$ \Box f = (Qf \wedge \Theta)/\Theta^{2} . $$
Then we have
\begin{lem}
$ \Box f= \Delta_{\Theta} f + V f $
where $V= h^{3} \Delta_{\omega}( h^{-1})$ and we are writing $\Delta_{\Theta}, \Delta_{\omega}$ for the Laplace operators of the two metrics.
\end{lem}

To see this, suppose $fmg$ have compact support and write,
$$  \int \Box f g d\mu= \int h \cD(h^{-1} f) \wedge g\Theta= \int \cD(h^{-1} f) (h^{-1} g) \omega. $$
Now apply Stokes' Theorem and the fact that $\omega$ is closed to write this as
$$  -2i \int \partial(h^{-f}) \db(h^{-1} g) \wedge \omega. $$
Some further manipulation, which we leave as an exercise for the reader, shows that this is equal to
$$  \int (\nabla f,\nabla g)+ V(fg) d\mu, $$
(where the inner product is computed using $\Theta$) with the stated function $V$. 

\

The conformal equivalence the  the flat metric $\Omega$ on $\bC^{2}\setminus\{0\}$ to the cylindrical metric corresponds to  $h=r$. We then have 
$$\Delta_{\Omega}r^{-1}= r^{-3} \frac{\partial}{\partial r}(r^{3} \frac{\partial r^{-1}}{\partial r})= r^{-3}, $$
so in this case $V=1$. 

\

Now return to our manifold $Z$ with the Kahler metric $\omega_{0}$ depending on the parameter $R$. We have a function $r_{Y}$ on each copy of $Y$. Let $r_{X}$ be a positive function on $X$ which, in a fixed ball about each singular point, is equal to the distance to that singular point. There is then a function
$h$ on $Z$ equal to $r_{X}$ on the \lq\lq $X$-side'' and to $R^{-1} r_{Y}$ on the \lq\lq $Y$-side''. (Since we glued the metrics on the sphere where
$r_{Y}=R^{1/2}, r_{X}=R^{-1/2}$.) The hermitian metric $\Theta_{0}=h^{-1} \omega_{0}$ contains a long cylinder. More precisely there is a region in $Z$ which we can identify with a cylinder $P^{3}\times (-T,T)$ where $P^{3}=S^{3}/\pm 1$ and $T$ is approximately $(\log R)/2$. On this cylinder the co-ordinate $t\in (-T,T)$ is $\log R^{1/2} h$. The metric $\Theta_{0}$ is precisely cylindrical on the part of the cylinder $t\geq 0$ and is approximately cylindrical on the region $t\leq 0$. We can think of $\Theta_{0}$ as being obtained in the following way. Define  the Hermitian metric $\Theta_{X}= r_{X}^{-1} \omega_{X}$ on $X$: this has 16 cylindrical ends. Now take the metric $\Theta_{Y}= r_{Y}^{-1}\omega_{Y}$ on $Y$. This is a metric with an asymptotically cylindrical end. Now \lq\lq cut-off'' the metric $\Theta_{Y}$ at a  distance $T/2$ down the end, to make it exactly cylindrical, and perform the \lq\lq connected sum'' construction considered before (except, of course, that we have 16 copies of $Y$).  This \lq\lq cutting off'' is exactly what we have specified above, but we are now viewing it in a slightly different way. We have a differential operator $\Box$ on $Z$ which is equal to $\Delta_{\Theta_{0}}+V$ where $V$ is equal to $1$ in the region $t>0$ of the cylinder can be supposed to be close to $1$ on the region $t<0$. Again, we have corresponding operators on the complete manifolds $X,Y$ with asymptotically tubular ends. The point of all this is that we can apply our general analytical theory to the operator $\Box$. (For this we need to extend the discussion, as mentioned before, to metrics with asymptotically cylindrical ends,  but that presents no difficulties)

\subsection{The proof}

Everything is now in place to proceed with the proof. We suppose our metric $\omega_{X}$ is chosen so that $\omega_{X}^{2}=\chi\wedge \overline{\chi}$. We seek a function
$\phi$ on $Z$ and $\lambda>0$ such that 
$$  (\omega_{0}+\cD \phi)^{2}= \lambda \chi\wedge \overline{\chi}$$
which is to say
$$  (\omega_{0}+ \cD \phi)^{2} = \lambda (1+\eta) \omega_{0}^{2}, $$
where we make a stretch of language to consider the function $\eta$ as a function on $Z$ in the obvious way. In the cylindrical picture $\eta$ is supported on a band $\vert t\vert \leq \log 2$ say (or, really, 16 such bands, one for each gluing region). We have $\vert \eta\vert =O(R^{-2})$ and it is easy to see that the same holds for all derivatives of $\eta$. So for any $k$ the $L^{2}_{k}$ norm of $\eta$ is $O(R^{-2})$. Now write
$\phi=hf$ and express the equation in terms of $\Theta_{0}= h^{-2} \omega_{0}$. We get
$$  (\Theta_{0} + h^{-2} \cD( h f))^{2}= \lambda (1+\eta) \Theta_{0}^{2}. $$
Expanding the quadratic term and rearranging, this is
$$    \Box f + h^{-3} Q(f)^{2} = h^{3} (\lambda(1+\eta)-1). $$
The problem here is that $h$ is very small on the \lq\lq $Y$-side'', in fact $O(R^{-1})$, so the co-efficient of $Q(f)^{2}$ is very large. To deal with this, set $f=R^{-3} g$.  So we have an equation for the pair $(g,\lambda)$ which is
$$  \Box g + (Rh)^{-3} Q(g)^{2} = (Rh)^{3} (\lambda(1+\eta)-1)  . $$
Now $(Rh)^{-1}$ is bounded (along with all its derivatives). The differential operator $Q$ has co-efficients which are bounded, along with all derivatives independent of $R$. We will solve the equation for  $g$ in the Sobolev space $L^{2}_{5}$. Then our multiplication $L^{2}_{3}\times L^{2}_{3}\rightarrow L^{2}_{3}$ implies that
 
$$  \Vert (Rh)^{-3} ( Q(g_{1})^{2} -Q(g_{2})^{2}) \Vert_{L^{2}_{3}} \leq C \Vert g_{1}- g_{2} \Vert_{L^{2}_{5}} \ \left( \Vert g_{1}\Vert_{L^{2}_{5}}+ \Vert g_{2}\Vert_{L^{2}_{5}}\right). $$

Also since $L^{2}_{3}\subset C^{0}$, a small solution $g$ in $L^{2}_{5}$ will define a positive form. So everything is in place to try to apply the inverse function theorem.

Now $\eta$ is supported  in a  band of fixed width in the middle of the cylinder and on this band $Rh$ is $O(R^{-1/2})$. Since $\eta$ is $O(R^{-2})$ we see that $(Rh)^{3}\eta$ is $O(R^{-1/2})<<1$. The same holds for all derivatives.
It is time to examine the linearised problem which, by our general theory, reduces to considering  the kernel of the operator
$\Box$ over the complete manifolds $X,Y$.

  \
  
  By the definition of $\Box$, a function $f$ satisfies $\Box f=0$ if and only if $\Delta_{\omega} (h^{-1} f)=0$. Consider first $f$ on $Y$. Then if $f$ tends to zero at infinity the same is true {\it a fortiori} for $r_{Y}^{-1} f$ and if $\Delta_{\omega} (r_{Y}^{-1} f)=0$ the function must vanish by the maximal principle. So there is no kernel of $\Box$ on $Y$. Similarly, a function in the kernel of $\Box$ on $X$ corresponds to a harmonic function, in the flat metric, which is $o(r_{X}^{-1})$ at the singularities. Since the fundamental solution of the Laplacian is 4 dimensions is $O(r^{-2})$ the only possibility is that this function is constant. So there is a 1-dimensional kernel of $\Box$ on $X$, spanned by the function $r_{X}$. Indeed there is obviously a kernel of $\Box$ on $Z$, spanned by the function $h$. Thus we are in a slighly more complicated situation than that envisaged before, but the same argument shows that we can invert $\Box$ on $Z$ uniformly \lq\lq modulo $h$''. That is there is a uniformly bounded operator $P$ and a linear functional $\pi$ such that
$$    \rho= \Box P( \rho) + \pi(\rho) h. $$

This fits in with the fact that we have an additional parameter $\lambda$ in our problem. Going back to the Kahler picture we know that the metrics $\omega_{0}$ and $\omega_{0}+\cD\phi$ on $Z$ have the same volume. This goes over to the identity
$$   \int h( \Box g + (Rh)^{-3} Q(g)^{2}) d\mu =0, $$
for any $g$. Thus the parameter $\lambda$ is determined by $\eta$ through the equation
$$ \lambda \int_{Z} (1+\eta)h^{4}\  d\mu = \int_{Z} h^{4}\ d\mu. $$
We define $\lambda$ by this formula, so $\lambda=1+O(R^{-4})$, since $h=O(R^{-1/2})$ on the support of $\eta$. With this value of $\lambda$ we  solve the nonlinear equation \lq\lq modulo $h$'' by the inverse function theorem. That is, we  solve the equation for $(g,\tau)$, where $\tau$ is a constant,
$$\Box g + (Rh)^{-3} Q(g)^{2} = (Rh)^{3} (\lambda(1+\eta)-1) + \tau h. $$
Now taking the $L^{2}$ inner product with $h$ we see that in fact $\tau=0$ and we have found our Calabi-Yau metric.

%+Bibliography

%-Bibliography

\end{document}